\documentclass{article}
\usepackage{amssymb,graphics}
\newcommand{\bbZ}{{\mathbb Z}}
\newcommand{\bbQ}{{\mathbb Q}}

\newcommand{\qed}{\mbox{$QED$}\vspace{\baselineskip}}
\newtheorem{cor}{Corollary}[section]
\newtheorem{lem}[cor]{Lemma}
\newtheorem{thm}[cor]{Theorem}

\newtheorem{prop}[cor]{Proposition}

\begin{document}

\title{Counting faces of cubical spheres modulo two}
\author{Eric K. Babson\thanks{Supported in part by an NSF
Postdoctoral Fellowship}
\ and Clara S. Chan\thanks{Supported in part by an NSF Career Grant}\\
Cornell University, Ithaca, NY 14853-7901\\
Wesleyan University, Middletown, CT 60459-0128}
\maketitle

\begin{abstract}
Several recent papers have addressed the problem of characterizing the 
$f$-vectors of cubical polytopes.  This is largely motivated by the complete 
characterization of the $f$-vectors of simplicial polytopes given by Stanley, 
Billera, and Lee 
in 1980 \cite{S} \cite{BL}.  Along these lines Blind and Blind \cite{BB} have shown that 
unlike in the 
simplicial case, there are parity restrictions on the $f$-vectors of cubical 
polytopes.  In particular, except for polygons, all even dimensional 
cubical polytopes must have an even number of vertices.  Here this result is 
extended to a class of zonotopal complexes which includes simply connected 
odd dimensional manifolds.  This paper then shows that the only 
modular equations which hold for the $f$-vectors of all d-dimensional 
cubical polytopes (and hence spheres) are modulo two.  Finally, the question of 
which mod two 
equations hold for the $f$-vectors of PL cubical spheres 
is reduced to a question about the Euler characteristics 
of multiple point loci from codimension one PL immersions into the $d$-sphere.
Some results about this topological question are known \cite{E} \cite{H} \cite{L} 
and Herbert's result we translate into the cubical setting, thereby 
removing the PL requirement.   A central definition in this paper is
that of the derivative complex, which
captures the correspondence between cubical spheres and codimension one
immersions.
\end{abstract}

We would like to thank Louis Billera for his contributions to this paper 
and Louis Funar for pointing us to Lannes' paper. 

\begin{section}{Introduction}
Several recent papers have addressed the problem of characterizing the 
$f$-vectors of cubical polytopes \cite{A} \cite{BBC} \cite{BCL} \cite{BB}
\cite{He} \cite{J}.
A {\em cubical $d$-polytope} is a $d$-dimensional
convex polytope all of whose boundary faces are combinatorially equivalent
to cubes. The $f$-vector of a cubical complex $K$ is the vector $f(K) = 
(f_0,f_1,f_2,\ldots,f_{d-1})$ or $ \sum f_it^i$, where $f_i$ denotes the number of
$i$-dimensional faces in $K$. The motivation for characterizing the
$f$-vectors of cubical polytopes stems from the success of
the characterization effort for $f$-vectors of simplicial polytopes \cite{S} \cite{BL}.

One of the first cubical $f$-vector results was that of Blind and Blind \cite{BB},
who showed that unlike in the simplicial case, there are parity restrictions
on the $f$-vectors of cubical polytopes.  In particular, except for polygons,
all even dimensional cubical polytopes must have an even number of vertices.  
In Section 2 we extend this result to a class of complexes including all
simply connected odd dimensional zonotopal manifolds.  The case of particular
interest to us is cubical spheres.

In Section 3 we show that 
the only modular equations which hold for the $f$-vectors of all
d-dimensional cubical polytopes are modulo two.  Since the linear equations 
which hold for spheres and polytopes are the same, this result also holds 
for cubical spheres.

In the remainder of the paper we focus on piecewise linear (PL) cubical
spheres, which are cubical complexes whose realizations are PL-homeomorphic 
to the boundary of a simplex \cite{RS}.  This is a more general class of complexes 
than boundaries of cubical polytopes.

In Section 4, we define the derivative complex of any cubical complex.  
For a PL cubical $d$-sphere $K$ the derivative complex is a codimension one 
PL cubical manifold with an immersion into $|K|$ such that the image is the 
$(d-1)$-skeleton of the dual to the cubical structure on $K$.  
This construction captures the correspondence
between cubical spheres and normal crossing codimension one immersions which
is crucial to the remainder of the paper.  It is called the derivative
complex because it acts as a derivation with respect to products and 
disjoint union.

In Section 5, we use the concept of derivative complexes to prove a new modulo
two equation for $f$-vectors of cubical spheres analogous to a topological
result about the corresponding immersions \cite{E} \cite{H}.

In Section 6, we construct a non-canonical PL cubical sphere from any 
codimension one normal crossing PL immersion into the sphere.  These are 
related by the fact that modulo two the number of $k$-cubes in the complex 
is the Euler characteristic of the $k$-fold self-intersection set of the 
immersion.  Combining this with the derivative complex of Section 4 
shows that the set of $f$-vectors of PL cubical $d$-spheres is the same 
modulo two as the set of Euler characteristics of multiple point
loci from normal crossing codimension one PL immersions into the $d$-sphere.

In Section 7 we use this correspondence to translate other results about 
immersions to results about piecewise linear cubical spheres. For these results 
we do not know if the PL requirement can be removed.

Throughout the paper we will use standard poset terminology \cite{S2} 
chapter 3, which we review here.  

For a poset $P$ and $x \in P$, we denote by
$\bigwedge x$ the  {\em order ideal} $\{z\in P : z\le x \}$,
and by $\bigvee x$ the {\em filter} $\{z\in P : z\ge x \}$.  
We define the {\em link\/}\ of $x$ in $P$ to be the
poset $\bigvee x \backslash \{ x \}$.
$P^{op}$ denotes the dual poset to $P$, {\em i.e.,\/}\ the underlying
set of $P$ with the reverse order.
We denote by $|P|$ the (simplicial) complex of chains in $P$.
A map between posets is called a poset map if it preserves order.

A poset is {\em cubical\/}\ if each order ideal $\bigwedge x$
is a product of copies of $I$, the three element face poset of an 
interval, excluding the empty set (so $I$ has a maximum but no minimum
element).

Cubical posets are ranked, the rank of an element being the number of links
in a maximal chain ending at this element.  Thus in a cubical complex,
rank is the same as dimension.  A cubical complex is a cell complex whose
face poset is a cubical poset $P$
such that $\widehat P$ ({\em i.e.,\/}\ $P$ with a minimum and maximum element
adjoined) is a lattice.
Throughout this paper when we consider the face poset of any
cubical complex we will always exclude the empty set.

To keep notation to a minimum, we will usually denote a cell complex
and its face poset by the same symbol.
An element of a face poset is called an $i$-face, if it has rank $i$.
If two faces $F$ and $G$ are related by $F< G$ then $F$ is called a
face of $G$, and if also $rank (F) = rank (G) - 1$ then $F$ is called
a facet of G.
A flag of faces is a set of faces which are totally ordered, {\em i.e.,\/}\
$F_1 < F_2 < \ldots < F_k$.
The Euler characteristic $\chi (P)$ of a face poset $P$ is the alternating sum
of the number of faces of each rank.

\end{section}

\begin{section}{Bicolorings and Even Vertices}  
We begin by generalizing the result of Blind and Blind which says that apart 
from polygons, every even dimensional cubical polytope must have an even 
number of vertices.  

The proof entails first showing that apart from circles, cubical spheres 
are bicolorable and then that every bicolorable odd dimensional sphere 
has the same number of vertices of each color.  Here a bicoloring of a 
complex is a choice from two colors (e.g. black and white) for each 
vertex so that each edge has one vertex of each color, ({\em i.e.} the 1-skeleton 
is bipartite).

\begin{prop}
If $K$ is a cubical complex with $H^1(K;\bbZ/2\bbZ)=0$, then $K$ has a 
bicoloring.  
\end{prop}

{\em Proof:\/}\
Note that if $K^1 \in C^1(K;\bbZ/2\bbZ)$ is the cellular cochain which assigns 
one to every
edge, then a bicoloring of $K$ is a cochain $c \in C^0(K;\bbZ/2\bbZ)$ such that  
$\delta c=K^1$.
Since every 2-face of $K$ has an even number of edges (four),
$K^1$ is a cocycle.  Thus $K$ has a bicoloring if and
only if $[K^1]$ is 0 in $H^1(K;\bbZ/2\bbZ)$.  In particular, if 
$H^1(K;\bbZ/2\bbZ)$ is trivial, 
then $K$ must have a bicoloring.  
\qed

  A complex is $n$-Eulerian if the link of every vertex has 
Euler characteristic $n$.  

\begin{thm}
If $K$ is a bicolored $n$-Eulerian cubical complex with $n$ nonzero, then 
the same number of vertices are assigned each color, and hence $K$ has an
even number of vertices.
\end{thm}

{\em Proof:\/}\
Let $f_b$ and $f_w$ denote the number of vertices of $K$ colored black and 
white, respectively.  Let $f_{b,i}$ and $f_{w,i}$ denote the number of 
pairs $(v,k)$ where $v$ is a vertex of the $i$-face $k\in K$ and $v$ is colored 
black or white, respectively.  
Since $K$ is $n$-Eulerian, the sum over the links of all black vertices gives 
$$ n \cdot f_b=f_{b,1}-f_{b,2}+ \ldots +(-1)^{d-1}f_{b,d}. $$ 
Similarly for the white vertices.  
If a cube  has dimension $i\ge 1$, the bicoloring of its vertices has 
the same number of each color, so we get $f_{b,i}=f_{w,i}$ for every 
$i \ge 1$.  
Thus the two equations above give $n \cdot f_b=n \cdot f_w$ and, since $n$ 
is nonzero, this gives $f_b=f_w$.  
\qed

\begin{cor}
If $K$ is a cubical $d$-sphere with $3 \leq d$ odd then $K$ has an even number 
of vertices.  
\end{cor}

Note that the proofs of the above proposition and theorem hold equally well 
if cubical complexes are replaced with zonotopal complexes.  Thus every 
$n$-Eulerian zonotopal complex with $n \neq 0$ and trivial $\bbZ/2\bbZ$ first 
cohomology has an even number of vertices.  

\end{section}

\begin{section}{Modulo two}

This section shows that the parity requirement on $f_0$ for odd-dimensional
cubical polytopes is not unusual, in fact every modular equation which holds 
for all $f$-vectors of cubical $(d+1)$-polytopes (or all cubical $d$-spheres) 
is modulo two.  
This restriction on the modulus follows from the main result of this section,
that the $\bbZ$-affine span $L^d$ of all $f$-vectors of boundaries of cubical 
$(d+1)$-polytopes contains a full rank affine sublattice $R^d$ consisting of 
all vectors in the $\bbQ$-affine span of $L^d$ with only even entries.  
We show that $R^d$ is in fact generated by the $f$-vectors of the boundaries 
of $(d+1)$-dimensional cubical zonotopes.

A { \em zonotope} is a polytope which can be generated by taking the 
Minkowski sum of a finite set of line segments \cite{Z} chapter 7.
If a minimal generating set for a zonotope contains $n$ line segments
then the zonotope is said to have $n$ {\em zones.}  A { \em cubical zonotope} 
is a zonotope all of whose boundary faces are combinatorially equivalent to cubes.
Let $F^d_n$ denote the $f$-vector of the boundary of any cubical
$(d+1)$-zonotope with $n+d+1$ zones, and $F^0_n=2$.  
Note that there are many such cubical zonotopes but all have the same $f$-vector.  
(This follows easily from the correspondence between cubical zonotopes
and generic hyperplane arrangements \cite{Za}.)
In particular $F^d_0=(2+t)^{d+1}-t^{d+1}$ is the $f$-vector of the boundary of the 
$(d+1)$-cube.

Since $F^d_0 \in R^d \cap L^d$, we can form the lattices
$R^d-F^d_0 = \{ y \vert y+F^d_0\in R^d \}$ and
$L^d - F^d_0 =   \{ y \vert y+F^d_0\in L^d \}$.

\begin{thm}
For every $d$, the lattice $R^d-F^d_0$ is a full rank sublattice of $L^d-F^d_0$.
\end{thm}

{\em Proof:\/}\
Let $E^d_n$ denote the $\bbZ$-affine span of $\{F^d_i\}_{i\ge n}$.

Note that $F^d_n=F^d_{n-1}+(1+t)F^{d-1}_n$ for each $n$.
Hence $E^d_n=\langle (1+t)E^{d-1}_{n+1} \rangle_{\bbZ} +F^d_n$, where
$\langle \ \ \rangle_{\bbZ}$ denotes $\bbZ$-linear span, and 
$\langle \ \ \rangle_{\bbQ}$ denotes $\bbQ$-linear span.  Thus, by induction on 
$d$, 
$E^d_n$ is independent of $n$, so write $E^d=E^d_n$ and we have
$E^d=\langle (1+t)E^{d-1} \rangle_{\bbZ} +F^d_0 $. 
To begin the induction note that $E^0_n=\{2\}$ for all $n$.  

We will show that rank $E^d = $ rank $L^d$, so $E^d-F^d_0$ is a full rank
sublattice of $L^d-F^d_0$.  Then we will show that $R^d = E^d$, completing the
proof.

Let $v:L^d \rightarrow \bbZ$ denote evaluation at -1. Then $v(F^d_0)=1+(-1)^{d}$ 
while $v(p)=0$ for any $p \in \langle (1+t)E^{d-1} \rangle_{\bbZ}$.   Thus if $d$ 
is even, the rank of $E^{d+1}$ equals the rank of 
$\langle E^d \rangle_{\bbZ}$ which is one more than the rank of 
$E^d$.  Since rank $E^0 = 0$ and rank $E^d \le $ rank 
$L^d = \lfloor \frac{d+1}{2} \rfloor$ for all $d\ge 1$ \cite{G}, we have rank 
$E^d = $ rank $L^d$ for all $d\ge 1$.

Now we show that $E^d=R^d$.  Since zonotopes are centrally symmetric,  
$E^d \subseteq R^d$, so we just need to show that $R^d\subseteq E^d$.
 It is clear that $R^1 \subseteq E^1$, so we proceed to
prove that $R^d \subseteq E^d$ by induction on $d$. 
Fix $w \in R^d - F^d_0$.  
Since $E^d-F^d_0$ is a full rank sublattice of $L^d-F^d_0$,  we have
$R^d-F^d_0 \subseteq \langle L^d-F^d_0 \rangle_{\bbQ}= \langle E^d-F^d_0 
\rangle_{\bbQ}=\langle (1+t)\langle E^{d-1}\rangle_{\bbZ}\rangle_{\bbQ}=(1+t) 
\langle E^{d-1} \rangle_{\bbQ}$.  Thus we 
can write $w\in R^d-F^d_0$ as $w = (1+t)\sum _i \lambda _i  w_i$, for 
some $w_i \in E^{d-1}$ and $\lambda _i \in \bbQ$.  
We will say that a vector is even if it has all even entries.
Since $w=(1+t)\sum_i \lambda _i w_i$ is even, 
$u =\sum_i \lambda _i w_i$ is also even.

If $d$ is even, then $\langle E^{d-1} \rangle_{\bbZ} = E^{d-1}$, so $u$ is an 
even
vector in the $\bbQ$-affine span of $E^{d-1}$.  In particular, this means that
$u$ is in $R^{d-1}$. But $R^{d-1} \subseteq E^{d-1}$ (by induction on $d$),
so we have
$u\in E^{d-1}$, which implies $w = (1+t)u \in E^d-F^d_0$.

If $d$ is odd, we use that $v(F^{d-1}_0)=2$ and $\langle E^{d-2} 
\rangle_{\bbZ}=E^{d-2}$ to show that 
$w\in E^d-F^d_0=\langle (1+t)(\langle (1+t)E^{d-2} \rangle_{\bbZ}+F^{d-1}_0) 
\rangle_{\bbZ}=(1+t)\langle(1+t)E^{d-2}+F^{d-1}_0 \rangle_{\bbZ}$.  
We can write $w=(1+t)u$ as above and $u = \mu _0 F^{d-1}_0 + (1+t)\sum _{i\ge 1} 
\mu _i v_i$ for some $v_i\in E^{d-2}$ and $\mu _i\in \bbQ$ with $\mu _0 = \sum 
_{i\ge 1}\mu _i$.
Since $u$ is even, $v(u)$ must be even, so 
$(1/2)v(u)=\mu _0 = \sum_{i\ge 1} \mu _i = n$ for some integer $n$.  
Now set $y =\sum_{i\ge 1} \mu_i v_i$ is an even vector in 
$\langle E^{d-2} \rangle_{\bbQ}$ and $\langle E^{d-2}\rangle_{\bbZ} = E^{d-2}$, 
so $y\in R^{d-2}$.  But $R^{d-2}\subseteq E^{d-2}$ (by induction on $d$),
so we have $y\in E^{d-2}$, which implies
$w =(1+t)(nF^{d-1}_0 + y) \in E^d-F^d_0$.

Thus we have that for any $w\in R^d-F^d_0$,
$w$ must be in $E^d-F^d_0$; so $R^d\subseteq E^d$, as desired.
\hfill
\qed

Note that the affine equations satisfied by all $f$-vectors of cubical 
polytopes are the same as those satisfied by all $f$-vectors of cubical 
spheres, since these equations can be derived from the Dehn-Sommerville
equations for simplicial spheres \cite{S1} by considering the links of 
faces in cubical polytopes and spheres.
Thus the $\bbZ$-affine span of all $f$-vectors of cubical $d$-spheres 
contains $L^d$ as a full rank affine sublattice, and so the theorem holds 
also for cubical $d$-spheres in place of the boundaries of cubical 
$(d+1)$-polytopes.  

\begin{cor}
The only modular equations which hold among the components of the
$f$-vectors of all cubical $d$-polytopes (and hence spheres) are modulo two.
\end{cor}

{\em Proof:\/}\
Let $V= \frac{1}{2}(R^d-F^d_0)$ denote the lattice of integral points in 
$\langle R^d - F^d_0 \rangle_{\bbQ}$.  
Then $L^d-F^d_0 \subset V$, since
$R^d-F^d_0$ is a full-rank sublattice of $L^d-F^d_0$.  A modulo $m$
equation on $L^d$ is a $\bbZ$-linear map $S:V\rightarrow \bbZ/(m\bbZ)$ which
restricts to 0 on $L^d-F^d_0$.  We may assume that $S$ is surjective,
otherwise we can get an equivalent surjective modular equation 
with modulus the order of the image of $S$.   
Thus $S(v)\equiv _m 1$ for some integral $v$ in 
$\langle \frac{1}{2}(R^d - F^d_0) \rangle_{\bbQ}$. 
Then $2v$ has all even entries and is in the same $\bbQ$-linear span,
so $2v\in R^d - F^d_0 \subset L^d - F^d_0$.   But then $2 \equiv _m S(2v) \equiv 
_m 0$, so $m=2$.  \qed

Note that a slight extension of the above argument shows
there are no modular restrictions on the $f$-vectors of arbitrary cubical 
$d$-manifolds -- since zonotopes are centrally symmetric, 1/2 
the $f$-vector of the boundary of any cubical zonotope is the $f$-vector 
of a cubical subdivision of real projective space.  

\end{section}

\begin{section}{Derivative Complexes}
In this section we give the more straightforward direction for the equivalence 
between PL (piecewise linear) cubical $d$-spheres and codimension one PL normal 
crossing immersions into the $d$-sphere.  The topological objects involved are 
described in (\cite{J} Section 6), where they are attributed to MacPherson 
and Stanley.

If $K$ is a cubical poset, define a new cubical poset $NK$ with elements the
ordered pairs $(b,c)\in K\times K$ such that the join of $b$ and $c$ covers
both, while $b$ and $c$ have no meet.   Thus $b$ and $c$ are opposite facets of 
their join.  The partial order on $NK$ is the partial order on 
$K$ taken component-wise.  Let $ \epsilon :NK \rightarrow NK$ denote the 
involution $\epsilon (b,c)=(c,b)$.  Then the {\em derivative complex\/}\  
of $K$ is the quotient poset $DK = NK/ \epsilon $.  Note that $NK$ and $DK$ 
are both cubical posets and $NK$ is simply a double cover of $DK$.  An element 
$\{b,c\}=(b,c)/\epsilon \in DK$ corresponds to a slice through the interior of the 
join of $b$ and $c$, 
parallel to $b$ and $c$.  An element  
$(b,c)\in NK$ corresponds to the side of $\{b,c\} \in DK$ which faces $b$
in $K$.  See Figure 1 for geometric realizations of $NK$ and $DK$.  

%K,NK, and DK for K = 2-cube, boundary of 3-cube, and 3-cube

\begin{center}
{\scalebox{.75}{\includegraphics{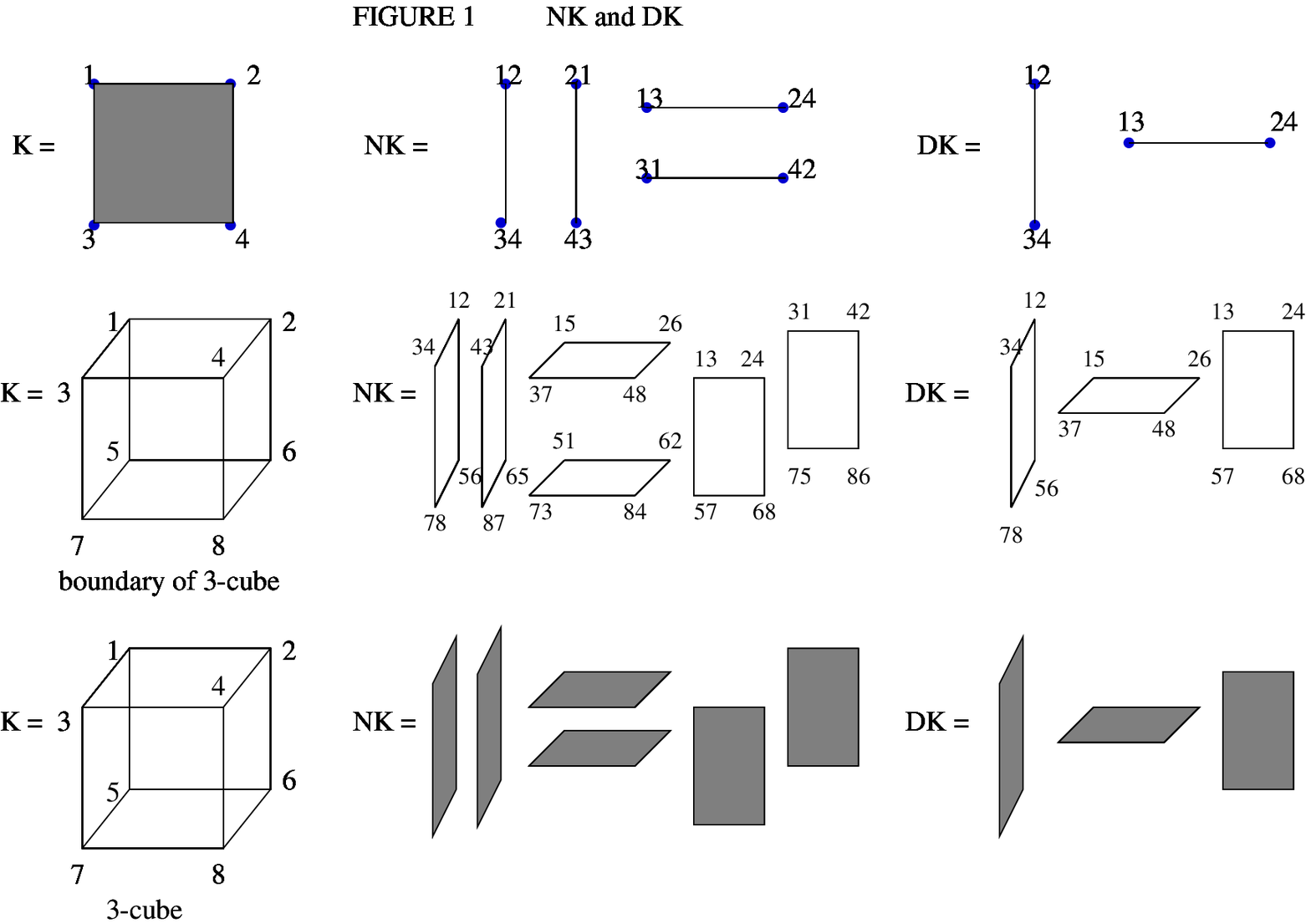}}}
\end{center}

Note that $$ f(DK,t) = \frac{d}{dt} f(K,t) $$
and $D$ and $N$ act as derivations with respect to product and disjoint union:
$$ D(K_1 \times K_2)=(DK_1 \times K_2) \cup (K_1 \times DK_2).$$

Finally we have the map $j : NK \rightarrow K$ 
taking $(b,c)$ to the join of $b$ and $c$, and similarly,
$j: DK \rightarrow K$.  
Both maps induce isomorphisms on the links of faces, so 
if $|K|$ is a piecewise linear manifold then so are $|DK|$ and $|NK|$.  
Furthermore $j ^*:|DK| \rightarrow |K|$ 
is a codimension one normal crossing immersion.
If dim$(|K|)=d$, then the image of $j^*$ is the $(d-1)$-skeleton of the dual to 
the cubical structure of $K$.
See Figure 2 for an illustration on the boundaries of cubical 3-polytopes.

%K and j(DK) for K = boundary of 3-cube and boundary of cubical octahedron

\begin{center}
{\scalebox{.75}{\includegraphics{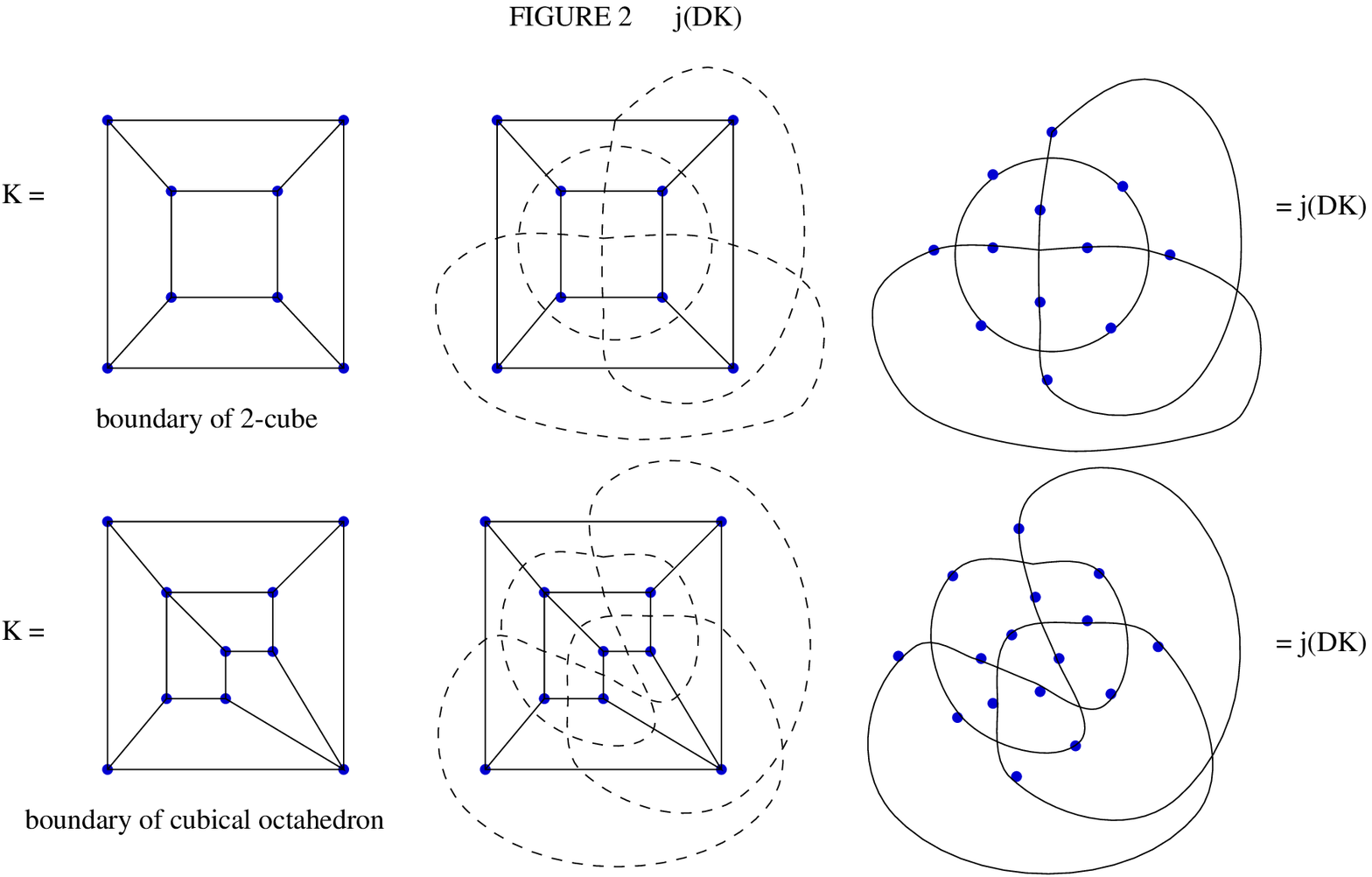}}}
\end{center}

In Section 6 we will show that any codimension one normal crossing 
immersion into $S^d$ is enumeratively equivalent, modulo two, to the immersion of 
a derivative complex.
Note that $f_i(K)=(-1)^{d-i} \chi (\{s\in |K|\vert |j^{*-1}(s)|=i\})$ since
the set of points $\{s\in|K|\vert |j^{*-1}(s)|=i\}$  is simply a disjoint
union of open $(d-i)$-balls, one for each $i$-face of $K$.  (Here
$\chi$ denotes Euler characteristic with closed supports, so $\chi
(R^n) = (-1)^n$.)
This relation will be used in Section 6.

\end{section}

\begin{section}{A Chain Lemma}

In this section we translate a known result about codimension one immersions
\cite{H} into a result about cubical complexes.  For even dimensional cubical 
polytopes of dimension at least 6, this result gives a modulo two condition 
different from the Blind-Blind \cite{BB} condition that $f_0$ must be even.  

If $K$ is a PL cubical $d$-sphere, then $j^*:|DK|\rightarrow |K|$ is
a codimension one normal crossing PL immersion into the $d$-sphere,
so a result of \cite{E}, \cite{H} tells us that if $d$ is odd, then
the number of degree $d$ intersection points is congruent modulo two to the Euler 
characteristic $\chi (|DK|)$.  In particular, $f_d(K) \equiv _2
f_0(DK) + f_1(DK) + f_2(DK) + \ldots f_{d-1}(DK)$.  Since $f_i(DK) = 
(i+1)f_{i+1}(K)$, this implies that $f_d(K) \equiv _2  f_1(K) + f_3(K) +
\ldots + f_d(K)$, or $f_1(K) + f_3(K) + \ldots f_{d-2}(K) \equiv _2 0$,
for all PL cubical $d$-spheres $K$ with $d$ odd.   Here we prove this
result directly for a class of cubical complexes including all
odd-dimensional cubical spheres, thus removing the PL requirement.

We first prove a lemma about simplicial flags of $K$, $DK$ and $NK$.
Some notation will be useful.  If $K$ is a poset, let $SK$ denote the simplicial 
poset of flags of $K$, and $C_*(K)=C_*(|K|; \bbZ/2\bbZ)$ the simplicial chain 
complex with 
coefficients in $\bbZ/2\bbZ$.  We identify elements of $C_*(K)$ with 
subsets of the set of elements of $SK$ so that both set and chain complex 
notations make sense
({\em e.g.,\/}\  $\cap, +, \partial$).

Let $\epsilon _* :C_*(NK) \rightarrow C_*(NK)$ denote the extension of the 
involution $\epsilon :NK \rightarrow NK$, and 
$\sigma =1+\epsilon _*: C_*(NK) \rightarrow C_*(NK)$.  One can check that 
$\partial \sigma =\sigma \partial $ 
and $\sigma (C_*(NK))$ is a chain subcomplex isomorphic to $C_*(DK)$.

Define a degree 0 map $\gamma =j^{-1}:C_*(K) \rightarrow C_*(NK)$ by
 $\gamma (a_0< \ldots  <a_r) = \sum _I (e_0< \ldots <e_r)$ where
 $(e_0< \ldots <e_r)\in I$ if and only if $j (e_i)=a_i$ for all $i$.
The image of $\gamma$ is contained in $\sigma (C_*(NK))$.  

Define a degree $-1$ map $\tau :C_*(K) \rightarrow C_*(NK)$ by 
$\tau(a_0< \ldots <a_r) = \sum _H (e_1< \ldots <e_r)$ where
$ (e_1< \ldots <e_r)\in H $ if and only if for every $i$,
$e_i=(b_i,c_i)$ such that $b_i\ge a_0$ and $a_i$ is the join of
$b_i$ and $c_i$.

Now suppose $e_1 <\ldots <e_{k-1}<e_{k+1}<\ldots <e_r$ and
for each $i$, $e_i=(b_i,c_i)$ such that $b_i\ge a_0$ and $a_i$ is
the join of $b_i$ and $c_i$.  Then 
if $k > 1$ there is a unique $e_k$ such that
$(e_1<\ldots <e_{k-1}<e_k<e_{k+1}<\ldots <e_r)\in H$,
but if $k=1$ then there exists such an $e_k$
if and only if $a_1 \not\le b_2$.
Based on this observation one can verify that $\partial \tau=\tau \partial$,
and also that $\sigma \tau= \partial \gamma +\gamma \partial $.

To prove the following lemma we often need to fix a side for each face
of $SDK$.  Such a choice corresponds to a chain $W\in C_*(NK)$ with 
$\sigma W=SNK$.  (Recall, since we identify elements of $C_*(NK)$ with
subsets of the set of elements in $SNK$, here we mean that $\sigma W$
is the chain of faces of $SNK$ with coefficient one for every face.)
Fix such a chain $W$.
For every $a=(a_1< \ldots < a_r) \in SK$, define $p_a\in C_*(K)$ as follows.
If $a_1$ is minimal in $K$, let $p_a = 0$ ({\em i.e.,\/},
the empty chain, by our identification of subsets of $SK$ with
elements in $C_*(K)$).  If $a_1$ is not minimal,
let $p_a=\{(p_0<a_1< \ldots <a_r)\}$ with $p_0$ minimal and $W \cap \gamma a
= \tau p_a$.  Note that since $p_0$ is minimal, $\gamma p_a=0$.

Figure 3 shows $P$, $\tau P$, and $\gamma P$, where $K$ is the solid
3-cube and $P$ is the boundary of a glued pair of triangles in the
subdivided 3-cube $SK$.  For this example, let 
$W$ choose the side of each face in $SDK$ which faces the facet with
smaller ternary labeling.  ({\em i.e.,\/}\
if $(b,c)\in NK$ occurs in a flag in $W$ then $b$ is less than $c$
with respect to the ternary labeling of faces shown.)
Then $W \cap \tau P = \{(001,201)\} + \{(002,202)\}$ and
$\sigma (W \cap \tau P) = \{(001,201)\} + \{(201,001)\} + \{(002,202)\} +
\{(202,002)\}$.

%$P$,$\tau P$, and $\gamma P$,where $P$ is the
%boundary of a glued pair of triangles in the subdivided 3-cube.
\begin{center}
{\scalebox{.75}{\includegraphics{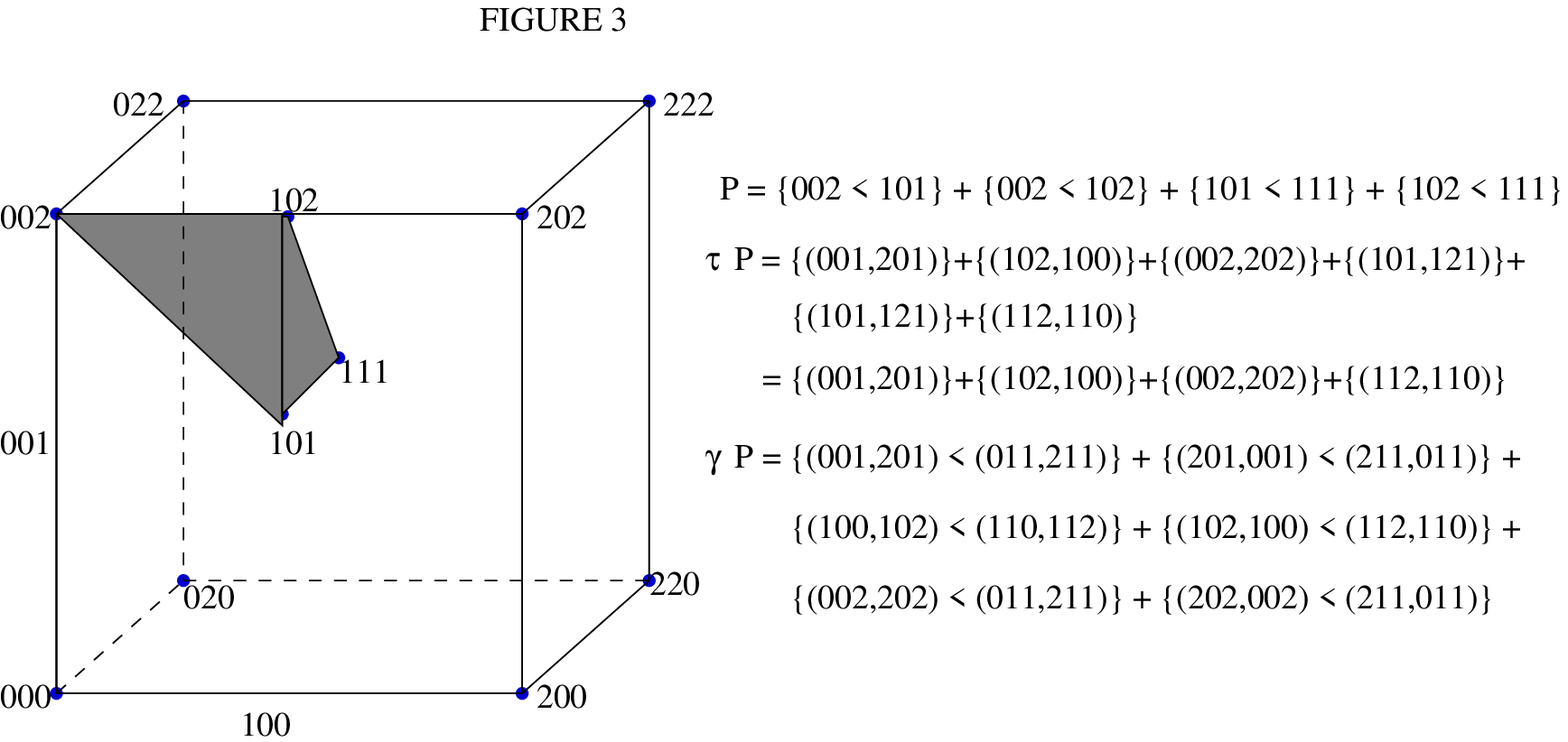}}}
\end{center}

\begin{lem}
If $K$ is a cubical complex and $P= \partial Q$ for some $Q \in C_*(K)$
and $\gamma P=\partial 
\sigma T$ for some $T\in C_*(NK)$ then $\tau P=\partial \sigma U$ for some $U\in 
C_*(NK)$.  
\end{lem}

{\em Proof:\/}\
Let $R = Q+ \sum _{a\in Q} \partial p_a$, so $\partial R = \partial Q = P$.
We also have $$\gamma R =
\gamma Q + \sum _{a\in Q} \gamma \partial p_a = \gamma Q + \gamma Q = 0,$$
since
$$\gamma \partial p_a=\partial \gamma p_a+ \sigma \tau p_a= \sigma \tau p_a =
\sigma (W \cap \gamma a)=\gamma a$$ and the coefficients are in $\bbZ/2\bbZ$. 
Thus $\sigma \tau R = (\gamma \partial +\partial \gamma ) R=\gamma P$.

Now note that for any $U \in C_*(NK)$, we can write
$U = \sigma A + B\cap \sigma U$ for some $A,B\in C_*(NK)$ with $\sigma B=SNK$.
So with $U=\tau R$, we have for some $A, B \in C_*(NK)$ with $\sigma B=SNK$, and 
$$\tau R = \sigma A + B \cap \sigma \tau R = \sigma A + B \cap \gamma P.$$
Thus
$$\tau P=\tau \partial R =\partial \tau R =
\partial (\sigma A + B\cap \gamma P)
 = \partial (\sigma A + B \cap \partial \sigma T)
 = \partial ( \sigma A + \sum _{t \in T} (B \cap \partial \sigma t)).$$ 

Now for each $t \in T$ there exist $A_t, B_t\in C_*(NK)$ with
$\sigma B_t = SNK$ and $\partial t\subseteq B_t$ and $B =B_t+\sigma A_t$.
So
$$B \cap \partial \sigma t=(B_t+\sigma A_t) \cap \partial \sigma t=
(B_t \cap \partial \sigma t) + \sigma (A_t \cap \partial t) = 
\partial t + \sigma (A_t \cap \partial t).$$
Thus
$$\tau P =
\partial(\sigma A + \sum_{t\in T}(\partial t + \sigma (A_t \cap \partial t)))
 = \partial \sigma U,$$
where $U = A + \sum_{t\in T}(A_t \cap \partial t) \in C_*(NK)$.
\qed

Now some careful choices together with the above lemma give us the desired
result.

 \begin{thm}
 If $K$ is an odd dimensional cubical $d$-sphere 
 then $f_1(K)+ f_3(K)+ \ldots +f_{d-2}(K)$ is even.  
 \end{thm}
 
 {\em Proof:\/}\ 
 We will show that there exists a set of edges in $SDK$ with
 $f_1(K)+f_3(K)+\ldots +f_{d-2}(K)$ boundary vertices (modulo 2),
 hence $f_1(K)+f_3(K)+\ldots +f_{d-2}(K)$ must be even.
 
 We will define a set of simplicial chains $P_i\in C_i(K)$ which
 satisfy the conditions of Lemma 5.1 by
 induction on $i$, so that finally we may conclude that $\gamma P_0$
 is a symmetric boundary in $C_*(NK)$, and from this we get the
desired zero boundary in $C_*(DK)$.
 
 Let $[d,d-i]=\{d,d-1,\dots,d-i\}$ and let $K_{[d,d-i]}$ denote
the $[d,d-i]$ rank-selected subposet of $K$.
% Also, let $(SK)_i=SK \cap C_i(K)$.  
Then define $P_i = (SK + S(K_{[d,d-i]})) \cap C_i(K)$.
 
 One can check that $\gamma P_i=(SNK+S(NK)_{[d-1,d-i-1]})\cap C_i(NK)$ by
 noting that every $(i+1)$-chain $e_0 < e_1 < \ldots < e_i$ in $C_i(NK)$
 has a unique preimage $\gamma ^{-1}(e_0 < e_1 < \ldots < e_i) =
 (j(e_0) < j(e_1) < \ldots < j(e_i))$ in $C_i(K)$.    
 It is also not hard to check that $ \tau P_i=\gamma P_{i-1}$ for any
 cubical complex, by noting that for every chain $e_1 < e_2 < \ldots < e_i$
 in $C_i(NK)$ there is an odd number of $(i-1)$-chains in
 $\tau ^{-1}(e_1 < e_2 < \ldots < e_i)$, in particular all chains of
 the form $a_0 < j(e_1) < j(e_2) < \ldots < j(e_i)$ where $a_0\le b_1$
 and $e_1=(b_1,c_1)$.  (Note that there is an odd number of such $a_0$
 since the number of faces in a cube, including the cube itself, is $3^t$
 where $t$ is the dimension of the cube.)  Since all coefficients are mod two,
 any odd number is equivalent to one.   
 This shows that $\tau P_i = (SNK+S(NK)_{[d-1,d-i]})\cap C_{i-1}(NK)$, so
 $\tau P_i = \gamma P_{i-1}$, as desired.
 
 Similarly, it is easy to check that $\partial P_i = 0$, since the
 preimage (under $\partial$) of any chain in $C_{i-1}(K)$ contains
 an even number of $i$-chains.  More explicitly, the first term of
$P_i$  has trivial boundary because 
 the Euler characteristic of the link of any face in a cubical sphere
 of any dimension is even, and a cube of any dimension has an even
 number of non-maximal, non-empty faces, while to check that the 
second term has trivial boundary uses that every codimension one 
face is in two facets, and that every $d-i-1$ cube contains an even 
number of facets.
 Thus, if $0<i<d$, then $\partial P_i=0$ so $H_i(|K|;\bbZ/2\bbZ)=0$ implies that 
$P_i= \partial Q_i$
 for some $Q_i \in C_{i+1}(K)$.  
 Finally, simply plugging in $i=d-1$ we get
 $\gamma P_{d-1}=0 =\partial \sigma 0$, and so by induction and Lemma 5.1,
if $0 \leq i <d$ then 
 $\gamma P_i=\partial \sigma U_i$ for some $U_i \in C_{i+1}(NK)$.  
 Thus $\gamma P_0 /\epsilon_*$ is a boundary and hence 
$|\gamma P_0/\epsilon_*|=\frac{1}{2}(\sum _j 2jf_j(K) + 2df_d(K))$ is even.  
Thus $f_1(K)+f_3(K)+\ldots +f_{d-2}(K) \equiv _2
 \sum _j jf_j(K)+ df_d(K) \equiv_2 0$, as desired.
 \qed
 
 Note that the proof holds more generally for any $d$ and any cubical complex 
 $K$ with no faces of dimension more than $d$, the link of every face having 
 even Euler characteristic, and $H_i(|K|, \bbZ/2\bbZ)=0$ for $0 <i <d$.
 Thus, for all such $K$, we have $df_d(K) \equiv _2 \sum _i i f_i(K)$.

\end{section}

 \begin{section}{Immersions to Cubations}
 
In Section 2, we constructed a codimension one PL normal crossing
immersion into the $d$-sphere which we called the derivative
complex, from any PL cubical $d$-sphere.
We now complete this correspondence.  In particular, we prove the following: 

\begin{thm}
Given a codimension one normal crossing PL immersion $y: M\rightarrow S^d$, there 
exists
a PL cubical $d$-sphere $K$ such that, modulo two, the Euler characteristics for
the multiple point loci of $y$ are the same as for the immersion of the
derivative complex, $j : |DK|\rightarrow |K|$.
In particular, $\chi (\{ s \in S^d||y^{-1}(s)|=i\}) \equiv _ 2 f_i(K)$.
\end{thm}

As a result, the PL case of our Theorem 5.2 proves the result
which motivated it \cite{E} \cite{H}:

\begin{cor}
If $d$ is odd and $y: M\rightarrow S^d$ is a codimension one normal crossing PL 
immersion, then the number of degree $d$ intersection points is congruent 
modulo two to the Euler characteristic $\chi (M)$.
\end{cor}

In the next section we will give other results about normal crossing 
immersions which  correspond to counting facets of PL cubical spheres 
modulo two.
 
Our proof of Theorem 6.1 involves combinatorial constructions of posets 
with PL-equivalent realizations.  These equivalences rely on regularity of
the cell complexes involved, as in the following lemma.
 
 \begin{lem}
 If $Q$ is the face lattice of a regular cell complex and 
$ g :P \rightarrow Q$ is a poset map so that for every $q \in Q$ we have that 
$|g^{-1}(\bigwedge q)|$ is a PL ball with dimension the rank of $q$ and 
boundary $|g^{-1}(\bigwedge q \backslash \{q\})|$ then $|P|$ is 
PL-equivalent to $|Q|$.  
\end{lem}
 
 {\em Proof:\/}\ 
 Construct a PL homeomorphism between $|P|$ and $|Q|$ by building 
 in one cell of $Q$ at a time, and completing the $k$-skeleton before 
 moving to $(k+1)$-cells.  
 To extend the map over a $k$-cell $q \in Q$ note that there is already a 
 chosen PL homeomorphism between $|\bigwedge q \backslash \{ q \}|$
 and $|g^{-1}(\bigwedge q \backslash \{ q \})|$, which are homeomorphic to    
 the $(k-1)$-sphere.
 We can extend this homeomorphism over $|g^{-1}(\bigwedge q)|$ to 
 $|\bigwedge q|$ by coning, as in \cite{RS} Lemma 1.10.
  \qed
 
{\em Proof of Theorem 6.1:\/}\ 
Given a PL codimension one normal crossing immersion $y:M \rightarrow S^d$, 
choose a triangulation $T$ of $S^d$ so that $M$ has a triangulation
making $y$ simplicial.  ($T$ exists by the definition of PL.)  
 
 The construction of the cubical sphere will be done in two steps, 
 each producing a PL sphere.  
 
Let $J$ be the poset with elements
$$\{(t,C_1,C)|t \in T, C_1 \subseteq C \in S(\bigwedge t)\}$$ 
and partial order given by $(t,C_1,C) \geq (t',C_1',C')$ if  
$t \geq t'$ and $C_1' \subseteq C_1 \subseteq C \subseteq C'$.  
Let $p:J\rightarrow T$  be the map $(t,C_1,C)\mapsto t$.
Note that the dual poset, $J^{op}$, is simply the cubical barycover \cite{BBC}
of the barycentric subdivision of $T$ with the subdivided faces
of the $(d-1)$-skeleton of $T$ thickened, as shown in Figure 4.
Since $T$ need not be dual cubical, $J^{op}$ need not be cubical.

%barycentric subdivision of a neighborhood of a 2-spx, cubical barycover 
%of that, with thickened (d-1)-skeleton of T, and dual of that.
\begin{center}
{\scalebox{.75}{\includegraphics{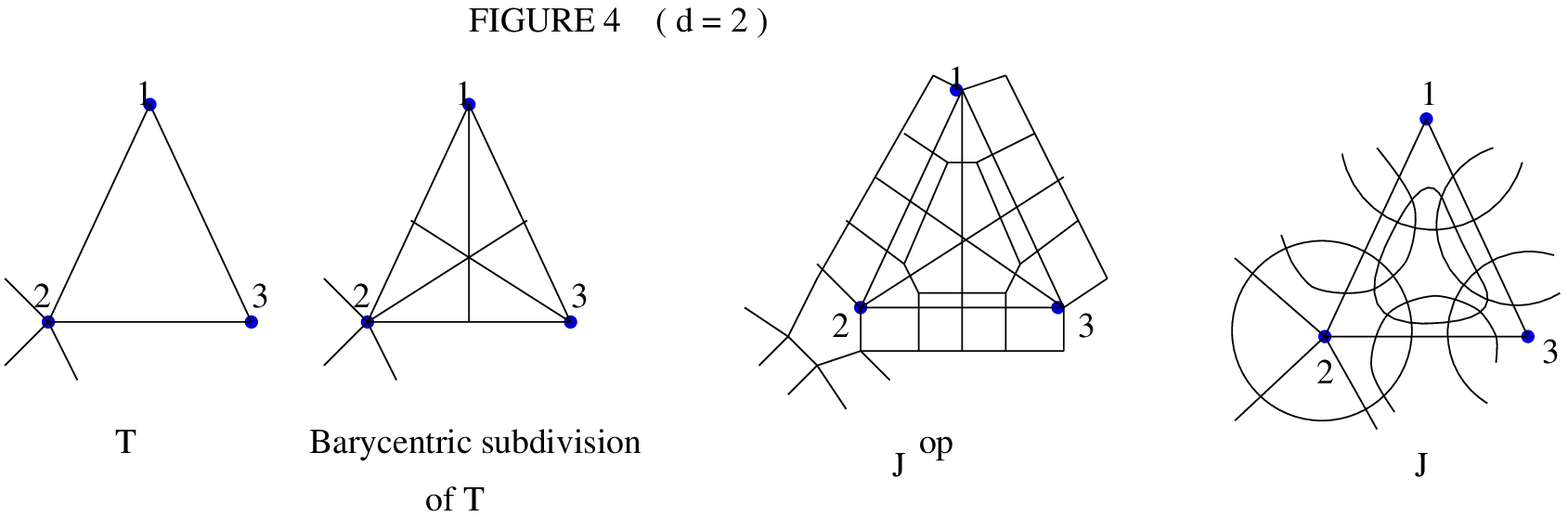}}}
\end{center}

One can check that $p:J\rightarrow T$ satisfies the conditions of the 
lemma.  Thus $|J|$ and $|T|$ are PL equivalent, so $|J|$ is a PL sphere.
Let $y:M\rightarrow |J|$ denote the induced simplicial
map of the associated subdivision of $M$ into $|J|$.  This again is a 
PL normal crossing immersion.  

Locally, $y(M)\subset |J|$ is a collection of normal crossing hyperplanes 
subdividing the open star of $t\in J$ into sectors, and hence there is a map 
$n_t$ taking each face in the open star of $t$ to its ambient sector.
Formally, for every $t \in J$ we define 
$n_t:\bigvee t \rightarrow (I^{op})^{|y^{-1}t|}$ where $I^{op}$ denotes the three 
element poset with a minimum but no 
maximum element.  Each copy of $I^{op}$ in the range of $n_t$ corresponds to
one of the normal crossing hyperplanes subdividing the open star of $t$,
in particular, the minimum element of $I^{op}$ corresponds to faces in that
hyperplane, and the other two elements of $I^{op}$ correspond to the two sides
of that hyperplane.
The maps $\{n_t\}$ are consistent, in the sense that if 
$ s \leq t$, then $n_s|_{\bigvee t}$ is constant and maximal on some copies 
of $I^{op}$ and agrees with $n_t$ on the other copies of $I^{op}$.

Let $K$ be the quotient of $J$ by the equivalence relation 
putting $(t,C_1,C) \sim (t',C_1,C)$ if $n_{MC}t=n_{MC}t'$,
where $MC$ is the maximum element of $C$.  
Then $K^{op}$ is cubical, as $K$ is the result of combining enough faces in $J$
that the only faces of $T$ which remain thickened are the normal crossing 
faces in $y(M)$.  See Figure 5.

%J restricted to a 2-simplex containing an edge of M <arrow> K restricted to
%the same 2-simplex.
\begin{center}
{\scalebox{.75}{\includegraphics{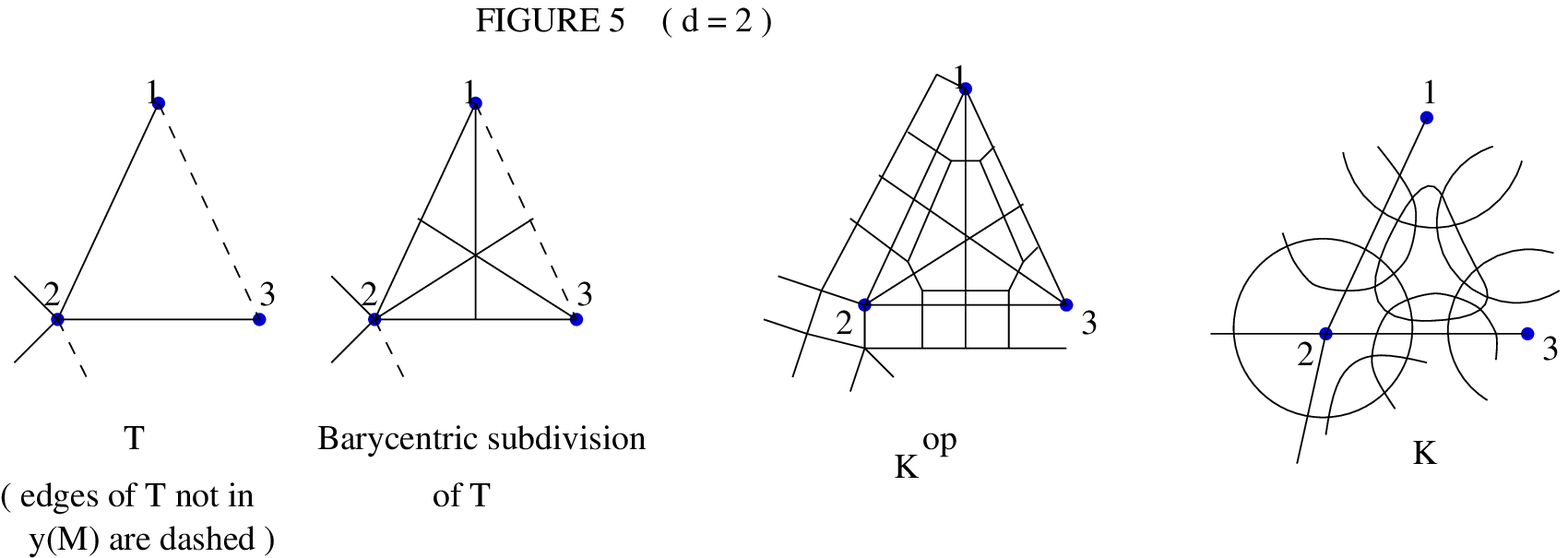}}}
\end{center}

Note that the fact that the immersion is normal crossing implies that 
$B=|\bigwedge n_t^{-1}r|$ is always a ball with boundary the union of 
$B \cap (\bigwedge (\bigvee t)\backslash \{\bigvee t\})$ and 
$n_t^{-1}( \bigwedge r \backslash \{r\})$.  
Thus the quotient map from $J$ to $K$ satisfies the conditions of the lemma,
so $K$ is a PL $d$-sphere and $K^{op}$ is a PL cubical $d$-sphere.

We now show that $f(K^{op})$ satisfies the desired modulo two equations.
For any $t\in T$, let $U(t)$ be the set of elements $(t',C_1,C)\in K$ for which $t 
= MC$, the maximum element of $C$.  See Figure 6.

%U(t) in a neighborhood of a 2-simplex.
\begin{center}
{\scalebox{.75}{\includegraphics{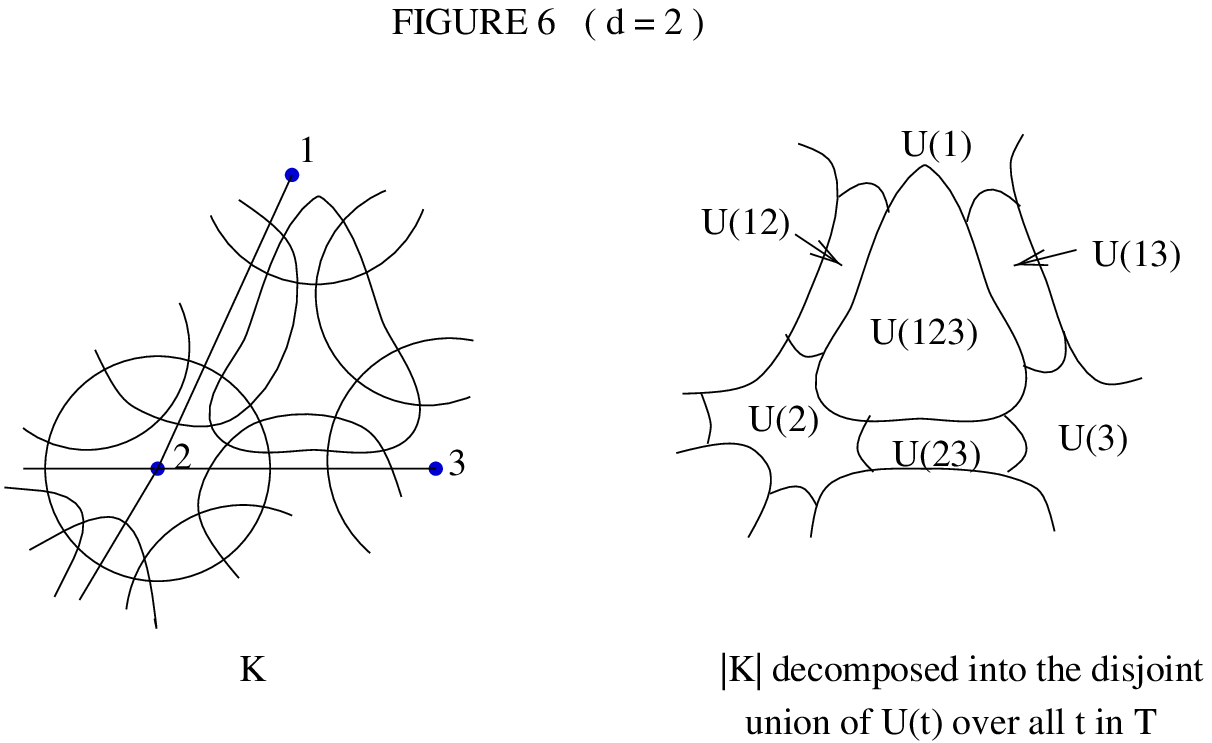}}}
\end{center}

Now note that there is an odd number of rank $n$ faces in $U(t)$ if and only 
if $n=d-|y^{-1}t|$.  
To see this, define an involution $i$ on $U(t)$ by 
$i(s,C_1,C) = (s', C_1',C')$ if the following hold:

1.  $C_1', C'$ are obtained from $C_1, C$ respectively by replacing each face $c$ 
strictly smaller than $t$ by its complement $t\backslash c$; and

2.  $n_t s$ agrees with $n_t s'$ on exactly those copies of $I^{op}$ for which 
either takes the minimum value, {\em i.e.,\/}\ $s$ and $s'$ are on opposite sides 
of all the local hyperplanes.
 
Then $i$ is rank preserving and fixes only the rank $d-|y^{-1}t|$ 
element $(t,\{t\},\{t\})$.

Thus $$\sum _{k \in K} x^{\hbox{rank}(k)}
              = \sum _{t \in T} \sum _{k \in U(t)}x^{\hbox{rank}(k)}
              \equiv _2\sum _{t \in T} x^{d-|y^{-1}t|},$$
and hence $$f_i(K^{op}) \equiv _2 \chi (\{ s \in S^d||y^{-1}(s)|=i\})$$ 
for all $i$, as desired.
\qed

\end{section}

\begin{section}{Facets and Questions}

Based on Theorem 6.1 and the smoothability of codimension one 
PL immersions (see \cite{HM} theorem 7.4) we can deduce the following about
the number of facets of PL cubical spheres.  
The first follows from \cite{B}, the second from \cite{L}, and the others 
from \cite{E}.  

1. There exists a PL cubical 3-sphere with an odd number of facets.  
(Let $j:RP^2\rightarrow S^3$ be Boy's immersion, having a single degree 3
intersection point \cite{B} \cite{E}.)  Thus the $\bbZ$-affine span of 
$f$-vectors of cubical 3-spheres is completely known, {\em i.e.,\/}
$f_0 \equiv _2 f_1 \equiv _2 f_2+f_3 \equiv _2 0$.

2. If $d$ is a multiple of 4 then a PL cubical $d$-sphere can have 
an odd number of facets if and only if $d=4$.  

3. If $d$ is odd, then a PL cubical $d$-sphere can have an
odd number of facets if and only if $d=1,3$ or $7$.  

4. If $d \equiv _4 2$, then a PL cubical $d$-sphere can have 
an odd number of facets only if $d=2^n - 2$ for some $n$.
Furthermore,  examples of cubical $d$-spheres with an odd number
of facets are known for $d=2,6,14,30,$ and $62$.  
(Eccles \cite{E2} showed that the existence of such spheres is equivalent to 
the existence of a framed $d$-manifold with Kervaire invariant one).

5. Edge orientable (in the sense of \cite{He}) PL cubical spheres 
can have an odd number of facets if and only if $d \in \{ 1,2,4\}$.\\

\noindent Some interesting open questions are:\\
Do the above hold for cubations of arbitrary topological spheres?\\  
Can even dimensional cubical spheres have an odd number of vertices?\\
For any given $d \ge 4$, what are all the modular equations
for $f$-vectors of cubical $d$-spheres?

\end{section}

\end{document}